\documentclass[10pt,twoside,reqno]{amsart}
\usepackage{amssymb}
\usepackage{multirow}
\calclayout

\numberwithin{equation}{section}
\makeatletter

\renewcommand{\@secnumfont}{\bfseries}

\renewcommand{\section}{\@startsection{section}{1}%
  {0mm}{.7\linespacing\@plus\linespacing}{.5\linespacing}
  {\normalfont\bfseries\centering}}

\newcommand{\bibsection}{\@startsection{section}{1}%
  {0mm}{.7\linespacing\@plus\linespacing}{.5\linespacing}
  {\normalfont\scshape\centering}}

\renewcommand{\@biblabel}[1]{#1.}

\newtheorem{thm}{\bf Theorem}[section]
\newtheorem{lem}[thm]{\bf Lemma}
\newtheorem{cor}[thm]{\bf Corollary}

\begin{document}

\vspace{1.3cm}

\title{On degenerate Central complete Bell polynomials}

\author{Taekyun Kim}
\address{Department of Mathematics, Kwangwoon University, Seoul 139-701, Republic
of Korea}
\email{tkkim@kw.ac.kr}

\author{Dae San Kim}
\address{Department of Mathematics, Sogang University, Seoul 121-742, Republic
of Korea}
\email{dskim@sogang.ac.kr}

\author{Gwan-Woo Jang}
\address{Department of Mathematics, Kwnagwoon University, Seoul 139-701, Republic of Korea}
\email{gwjang@kw.ac.kr}

\maketitle

\begin{abstract}
In this paper, we study the degenerate central complete and incomplete Bell polynomials which are degenerate versions of the recently introduced central complete and incomplete Bell polynomials and also central analogues for the degenerate complete and incomplete Bell polynomials. We investigate some properties and identities for these polynomials. 
\end{abstract}
\bigskip
\medskip

\section{\bf Introduction and preliminaries}
In this article, we consider the degenerate central incomplete Bell polynomials $T_{n,k}(x_1,x_2,\cdots,x_{n-k+1}|\lambda)$ given by
\begin{equation*}\begin{split}
&\frac{1}{k!}\Big(\sum_{m=1}^\infty x_{m}\big((\frac{1}{2})_{m,\lambda}-(-1)^{m}<\frac{1}{2}>_{m,\lambda}\big)\frac{t^m}{m!}\Big)^{k}\\
&\quad\quad\quad\quad\quad\quad=\sum_{n=k}^\infty T_{n,k}(x_{1},x_{2},\cdots,x_{n-k+1}|\lambda)\frac{t^n}{n!}, \,\,\,(\textnormal{see}\,\,\,\eqref{05},\eqref{06}),
\end{split}\end{equation*}
and the degenerate central complete Bell polynomials $B_n^{(c)}(x_1,x_2,\cdots,x_n|\lambda)$  given by
\begin{equation*} \begin{split} \label{27}
\textnormal{exp}\Big(&\sum_{i=1}^\infty x_{i}\big((\frac{1}{2})_{i,\lambda}-(-1)^{i}<\frac{1}{2}>_{i,\lambda}\big)\frac{t^i}{i!}\Big)\\
&\quad\quad\quad\quad\quad=\sum_{n=0}^\infty B_{n}^{(c)}(x_{1},x_{2},\cdots,x_{n}|\lambda)\frac{t^n}{n!}.
\end{split} \end{equation*}
and investigate some properties and identities for these polynomials. They are degenerate versions of the central complete and incomplete Bell polynomials. They are also viewed as 'central' analogues for degenerate complete and incomplete Bell polynomials (\textnormal{see}\,\,\,[9]), which are motivated by \eqref{04} and \eqref{10}. Before we move on to the next section, we will recall the necessary ingredients that are needed for our discussion in this paper. \\

For $n \geq 0$, the Stirling numbers of the first kind are given by
\begin{equation} \begin{split} \label{01}
\frac{1}{k!}\big(\log(1+t)\big)^{k} = \sum_{n=k}^\infty S_{1}(n,k)\frac{t^n}{n!},\,\,\,\,(\textnormal{see}\,\,\,[1-17]).
\end{split} \end{equation}

For $\lambda \in \mathbb{R}$, we define the degenerate exponential function as follows:
\begin{equation} \begin{split} \label{02}
e_{\lambda}^{x}(t)=(1+\lambda t)^{\frac{x}{\lambda}},\,\,\,\,(\textnormal{see}\,\,\,[2,14]).
\end{split} \end{equation}
Note that $\lim_{\lambda \rightarrow 0} e_{\lambda}^{x}(t)=\lim_{\lambda \rightarrow 0}(1+\lambda t)^{\frac{x}{\lambda}}=e^{xt}$. It is known that the degenerate Bell polynomials (also called degenerate Tochard polynomials or degenerate exponential polynomials) are defined by
\begin{equation} \begin{split} \label{03}
e^{x\big(e_{\lambda}(t)-1\big)}=\sum_{n=0}^\infty B_{n,\lambda}(x)\frac{t^n}{n!},\,\,\,\,(\textnormal{see}\,\,\,[15]),
\end{split} \end{equation}
where $e_{\lambda}(t)=e_{\lambda}^{1}(t)$.\\
When $x=1$, $B_{n,\lambda}=B_{n,\lambda}(1)$ are called the degenerate Bell numbers.\\
The degenerate incomplete Bell polynomials (also called degenerate partial Bell polynomials) are defined by the generating function (\textnormal{see}\,\,\,[9])
\begin{equation} \begin{split} \label{04}
\frac{1}{k!}\Big(\sum_{m=1}^\infty x_{m}(1)_{m,\lambda}\frac{t^m}{m!}\Big)^{k}=\sum_{n=k}^\infty B_{n,k}(x_{1},\cdots,x_{n-k+1}|\lambda)\frac{t^n}{n!},
\end{split} \end{equation}
where $k$ is a non-negative integer and $(x)_{m,\lambda}$ is the degenerate falling factorial sequence given by
\begin{equation} \begin{split} \label{05}
(x)_{0,\lambda}=1,\,\,\,\,(x)_{n,\lambda}=x(x-\lambda)\cdots\big(x-(n-1)\lambda\big),\,\,\,\,(n \geq 1).
\end{split} \end{equation}
Now, we define the degenerate rising factorial sequence as follows:
\begin{equation} \begin{split} \label{06}
<x>_{0,\lambda}=1,\,\,\,\,<x>_{n,\lambda}=x(x+\lambda)\cdots\big(x+(n-1)\lambda\big),\,\,\,\,(n \geq 1).
\end{split} \end{equation}
Note that $(-x)_{m,\lambda}=(-1)^{m}<x>_{m,\lambda}$.\\
From \eqref{04}, we note that
\begin{equation} \begin{split} \label{07}
B_{n,k}&(x_{1},x_{2},\cdots,x_{n-k+1}|\lambda)=\sum \frac{n!}{i_{1}!i_{2}!\cdots i_{n-k+1}!}\Big(\frac{x_{1}(1)_{1,\lambda}}{1!}\Big)^{i_{1}}\\
&\times \Big(\frac{x_{2}(1)_{2,\lambda}}{2!}\Big)^{i_{2}}\times \cdots\times\Big(\frac{x_{n-k+1}(1)_{n-k+1,\lambda}}{(n-k+1)!}\Big)^{i_{n-k+1}}, 
\end{split} \end{equation}
where the summation is over all integers $i_{1},i_{2},\cdots,i_{n-k+1}\geq 0$, such that $i_{1}+i_{2}+\cdots+i_{n-k+1}=k$ and $i_{1}+2i_{2}+\cdots+(n-k+1)i_{n-k+1}=n$.\\
It is known that the degenerate Stirling numbers of the second kind are defined by 
\begin{equation} \begin{split} \label{08}
\frac{1}{k!}\big(e_{\lambda}(t)-1\big)^{k}=\sum_{n=k}^\infty S_{2,\lambda}(n,k)\frac{t^n}{n!},\,\,\,\,(k \geq 0),\,\,\,\,(\textnormal{see}\,\,\,[8]).
\end{split} \end{equation}
From \eqref{04} and \eqref{08}, we note that
\begin{equation} \begin{split} \label{09}
B_{n,k}(\underbrace{1,1,\cdots,1}_{n-k+1-times}|\lambda)=S_{2,\lambda}(n,k),\,\,\,\,(n,k \geq 0).
\end{split} \end{equation}
By \eqref{07}, we easily get
\begin{equation*} \begin{split} 
B_{n,k}(\alpha x_{1},\alpha x_{2},\cdots,\alpha x_{n-k+1}|\lambda)=\alpha^{k}B_{n,k}(x_{1},x_{2},\cdots,x_{n-k+1}|\lambda),
\end{split} \end{equation*}
and
\begin{equation*} \begin{split} 
B_{n,k}(\alpha x_{1},\alpha^{2} x_{2},\cdots,\alpha^{n-k+1} x_{n-k+1}|\lambda)=\alpha^{n}B_{n,k}(x_{1},x_{2},\cdots,x_{n-k+1}|\lambda),
\end{split} \end{equation*}
where $\alpha \in \mathbb{R}$.\\
The degenerate complete Bell polynomials are defined by 
\begin{equation} \begin{split} \label{10}
\textnormal{exp}\Big(\sum_{i=1}^\infty x_{i}(1)_{i,\lambda}\frac{t^i}{i!}\Big)=\sum_{n=0}^\infty B_{n}(x_{1},x_{2},\cdots,x_{n}|\lambda)\frac{t^n}{n!}.
\end{split} \end{equation}
Then, by \eqref{04} and \eqref{10}, we get
\begin{equation} \begin{split} \label{11}
B_{n}(x_{1},x_{2},\cdots,x_{n}|\lambda)=\sum_{k=0}^{n} B_{n,k}(x_{1},x_{2},\cdots,x_{n-k+1}|\lambda).
\end{split} \end{equation}
From \eqref{11}, we note that
\begin{equation} \begin{split} \label{12}
B_{n}(x,x,\cdots,x|\lambda)&=\sum_{k=0}^{n} x^{k}B_{n,k}(1,1,\cdots,1|\lambda)\\
&=\sum_{k=0}^{n} x^{k} S_{2,\lambda}(n,k)=B_{n,\lambda}(x),\,\,\,\,(n \geq 0).
\end{split} \end{equation}
Recently, the degenerate central factorial numbers of the second kind are defined by 
\begin{equation} \begin{split} \label{13}
\frac{1}{k!}\big(e_{\lambda}^{\frac{1}{2}}(t)-e_{\lambda}^{-\frac{1}{2}}(t)\big)^{k}=\sum_{n=k}^\infty T_{2,\lambda}(n,k)\frac{t^n}{n!},\,\,\,\,(\textnormal{see}\,\,\,[11]),
\end{split} \end{equation}
where $k$ is a non-negative integer.\\
From \eqref{13}, we have
\begin{equation} \begin{split} \label{14}
T_{2,\lambda}(n,k)=\sum_{m=0}^{n}\Big(\frac{1}{k!}\sum_{l=0}^{k} {k \choose l}(-1)^{k-l}\big(l-\frac{k}{2}\big)^{m}\Big)\lambda^{n-m}S_{1}(n,m),
\end{split} \end{equation}
where $n,k \in \mathbb{Z}$ with $n \geq k \geq 0$, $(\textnormal{see}\,\,\,[11])$.\\
The degenerate central Bell polynomials are given by
\begin{equation} \begin{split} \label{15}
e^{x\big(e_{\lambda}^{\frac{1}{2}}(t)-e_{\lambda}^{-\frac{1}{2}}(t)\big)}=\sum_{n=0}^\infty B_{n,\lambda}^{(c)}(x)\frac{t^n}{n!}.
\end{split} \end{equation}
Thus, by \eqref{13} and \eqref{15}, we get
\begin{equation} \begin{split} \label{16}
B_{n,\lambda}^{(c)}(x)=\sum_{k=0}^{n} T_{2,\lambda}(n,k)x^{k},\,\,\,\,(n \geq 0).
\end{split} \end{equation}
When $x=1$, $B_{n,\lambda}^{(c)}=B_{n,\lambda}^{(c)}(1)$ are called the degenerate central Bell numbers.\\

\section{On degenerate central complete and incomplete Bell polynomials}
In view of \eqref{04}, we consider the degenerate central incomplete Bell polynomials given by 
\begin{equation} \begin{split} \label{17}
\frac{1}{k!}&\Big(\sum_{m=1}^\infty x_{m}\big((\frac{1}{2})_{m,\lambda}-(-1)^{m}<\frac{1}{2}>_{m,\lambda}\big)\frac{t^m}{m!}\Big)^{k}\\
&=\sum_{n=k}^\infty T_{n,k}(x_{1},x_{2},\cdots,x_{n-k+1}|\lambda)\frac{t^n}{n!},
\end{split} \end{equation}
where $k$ is a non-negative integer.\\
For $n,k \geq 0$ with $n-k \equiv 0$ $(mod\,\,\,2)$, by \eqref{17}, we get
\begin{equation} \begin{split} \label{18}
T_{n,k}&(x_{1},x_{2},\cdots,x_{n-k+1}|\lambda)=\sum \frac{n!}{i_{1}!i_{2}!\cdots i_{n-k+1}!}\\
&\times\Big(\frac{x_{1}\big((\frac{1}{2})_{1,\lambda}+<\frac{1}{2}>_{1,\lambda}\big)}{1!}\Big)^{i_{1}} \Big(\frac{x_{2}\big((\frac{1}{2})_{2,\lambda}-<\frac{1}{2}>_{2,\lambda}\big)}{2!}\Big)^{i_{2}}\times \cdots\\
&\times\Big(\frac{x_{n-k+1}\big((\frac{1}{2})_{n-k+1,\lambda}+<\frac{1}{2}>_{n-k+1,\lambda}\big)}{(n-k+1)!}\Big)^{i_{n-k+1}},
\end{split} \end{equation}
where the summation is over all integers $i_{1},i_{2},\cdots,i_{n-k+1}\geq 0$ such that $i_{1}+i_{2}+\cdots+i_{n-k+1}=k$ and $i_{1}+2i_{2}+\cdots+(n-k+1)i_{n-k+1}=n$.\\
From \eqref{18}, we can derive the following equation \eqref{19}.\\
For $n,k \geq 0$ with $n-k \equiv 0$ $(mod\,\,\,2)$, we have
\begin{equation} \begin{split} \label{19}
T_{n,k}&(x_{1},x_{2},\cdots,x_{n-k+1}|\lambda)\\
&=B_{n,k}\Big(x_{1}\big((\frac{1}{2})_{1,\lambda}+<\frac{1}{2}>_{1,\lambda}\big),x_{2}\big((\frac{1}{2})_{2,\lambda}-<\frac{1}{2}>_{2,\lambda}\big),\cdots,x_{n-k+1}\\
&\times\big((\frac{1}{2})_{n-k+1,\lambda}+<\frac{1}{2}>_{n-k+1,\lambda}\big)\Big).
\end{split} \end{equation}
Here $B_{n,k}(x_{1},x_{2},\cdots,x_{n-k+1})$ are the incomplete Bell polynomials which are defined by 
\begin{equation} \begin{split} \label{20}
B_{n,k}&(x_{1},x_{2},\cdots,x_{n-k+1})=\sum \frac{n!}{i_{1}!i_{2}!\cdots i_{n-k+1}!}\Big(\frac{x_{1}}{1!}\Big)^{i_{1}}\\
&\times \Big(\frac{x_{2}}{2!}\Big)^{i_{2}}\times \cdots\times\Big(\frac{x_{n-k+1}}{(n-k+1)!}\Big)^{i_{n-k+1}}, 
\end{split} \end{equation}
where the summation is over all integers $i_{1},i_{2},\cdots,i_{n-k+1}\geq 0$ such that $i_{1}+i_{2}+\cdots+i_{n-k+1}=k$ and $i_{1}+2i_{2}+\cdots+(n-k+1)i_{n-k+1}=n$.\\
Therefore, by \eqref{19} and \eqref{20}, we obtain the following lemma.
\begin{lem} For $n,k \geq 0$ with $n\geq k$ and $n-k \equiv 0$ $(mod\,\,\,2)$, we have
\begin{equation*} \begin{split}
T_{n,k}&(x_{1},x_{2},\cdots,x_{n-k+1}|\lambda)=B_{n,k}\Big(x_{1}\big((\frac{1}{2})_{1,\lambda}+<\frac{1}{2}>_{1,\lambda}\big),x_{2}\big((\frac{1}{2})_{2,\lambda}\\
&-<\frac{1}{2}>_{2,\lambda}\big),\cdots,x_{n-k+1}\big((\frac{1}{2})_{n-k+1,\lambda}+<\frac{1}{2}>_{n-k+1,\lambda}\big)\Big).
\end{split} \end{equation*}
\end{lem}
Let $n,k \geq 0$ with $n \geq k$ with $n-k \equiv 0$ $(mod\,\,\,2)$. Then, by \eqref{17}, we get
\begin{equation}\begin{split} \label{21}
\sum_{n=k}^\infty &T_{n,k}(x,x^{2},\cdots,x^{n-k+1}|\lambda)\frac{t^n}{n!}\\
&=\frac{1}{k!}\Big(x\big((\frac{1}{2})_{1,\lambda}+<\frac{1}{2}>_{1,\lambda}\big)t+x^2 \big((\frac{1}{2})_{2,\lambda}-<\frac{1}{2}>_{2,\lambda}\big)\frac{t^2}{2!}+\cdots\Big)^{k}\\
&=\frac{1}{k!}\big(e_{\lambda}^{\frac{1}{2}}(xt)-e_{\lambda}^{-\frac{1}{2}}(xt)\big)^{k}=\frac{1}{k!}e_{\lambda}^{-\frac{k}{2}}(xt)\big(e_{\lambda}(xt)-1\big)^{k}\\
&=\frac{1}{k!}\sum_{l=0}^{k} {k \choose l}(-1)^{k-l} e_{\lambda}^{(l-\frac{k}{2})}(xt)=\frac{1}{k!}\sum_{l=0}^{k} {k \choose l } (-1)^{k-l} \\
&\times e^{\frac{1}{\lambda}(l-\frac{k}{2})\log(1+\lambda xt)}
\end{split} \end{equation}
\begin{equation*} \begin{split}
&=\frac{1}{k!} \sum_{l=0}^{k} {k \choose l}(-1)^{k-l} \sum_{m=0}^\infty \lambda^{-m} \big(l - \frac{k}{2}\big)^{m} \frac{1}{m!}\big(\log(1+\lambda xt)\big)^{m}\\
&=\sum_{n=0}^\infty \bigg(\sum_{m=0}^{n} \Big(\frac{x^{n}}{k!} \sum_{l=0}^{k} {k \choose l}(-1)^{k-l} \big(l-\frac{k}{2}\big)^{m}\Big) \lambda^{n-m} S_{1}(n,m)\bigg)\frac{t^n}{n!}.
\end{split} \end{equation*}

Therefore, by comparing the coefficients on both sides of \eqref{21}, we obtain the following theorem.
\begin{thm} For $n,k \geq 0$ with $n-k \equiv 0$ $(mod \,\,\,2)$, we hvae
\begin{equation*} \begin{split} 
\sum_{m=0}^{n} \Big(\frac{x^{n}}{k!} &\sum_{l=0}^{k} {k \choose l}(-1)^{k-l} \big(l-\frac{k}{2}\big)^{m}\Big) \lambda^{n-m} S_{1}(n,m)\\
&=\begin{cases}T_{n,k}(x,x^2,\cdots,x^{n-k+1}|\lambda),\,\,\,\,{\rm if}\,\,\,\, n \geq k,\\
0,\,\,\,\,{\rm if}\,\,\,\,n < k.\end{cases}
\end{split} \end{equation*}
In particular,
\begin{equation*} \begin{split} 
\sum_{m=0}^{n} \Big(\frac{1}{k!} &\sum_{l=0}^{k} {k \choose l}(-1)^{k-l} \big(l-\frac{k}{2}\big)^{m}\Big) \lambda^{n-m} S_{1}(n,m)\\
&=\begin{cases}T_{n,k}(1,1,\cdots,1|\lambda),\,\,\,\,{\rm if}\,\,\,\, n \geq k,\\
0,\,\,\,\,{\rm if}\,\,\,\,n < k.\end{cases}
\end{split} \end{equation*}
\end{thm}
For $n,k \geq 0$ with $n-k \equiv 0$ $(mod\,\,\,2)$ and $n \geq k$, by \eqref{13} and \eqref{17}, we get
\begin{equation} \begin{split} \label{22}
T_{n,k}(1,1,\cdots,1|\lambda)=T_{2,\lambda}(n,k).
\end{split} \end{equation}
Therefore, by \eqref{22}, we obtain the following corollary.
\begin{cor} For $n,k \geq 0$ with $n-k \equiv 0$ $(mod\,\,\,2)$ and $n \geq k$, we have
\begin{equation*} \begin{split}
T_{n,k}(x,x^{2},\cdots,x^{n-k+1}|\lambda)=x^{n}T_{n,k}(1,1,\cdots,1|\lambda)=x^{n}T_{2,\lambda}(n,k),
\end{split} \end{equation*}
and
\begin{equation*} \begin{split} 
T_{2,\lambda}(n,k)&=T_{n,k}(1,1,\cdots,1|\lambda)=\sum \frac{n!}{i_{1}!i_{2}!\cdots i_{n-k+1}!}\\
&\times\Big(\frac{(\frac{1}{2})_{1,\lambda}+<\frac{1}{2}>_{1,\lambda}}{1!}\Big)^{i_{1}} \Big(\frac{(\frac{1}{2})_{2,\lambda}-<\frac{1}{2}>_{2,\lambda}}{2!}\Big)^{i_{2}}\times \cdots\\
&\times\Big(\frac{(\frac{1}{2})_{n-k+1,\lambda}+<\frac{1}{2}>_{n-k+1,\lambda}}{(n-k+1)!}\Big)^{i_{n-k+1}},
\end{split} \end{equation*}
\noindent where the summation is over all integers $i_{1},i_{2},\cdots,i_{n-k+1}\geq 0$ such that $i_{1}+i_{2}+\cdots+i_{n-k+1}=k$ and $i_{1}+2i_{2}+\cdots+(n-k+1)i_{n-k+1}=n$.\\
\end{cor}

For $n,k \geq 0$ with $n \geq k$ and $n-k \equiv 0$ $(mod\,\,\,2)$, we note from \eqref{21} that
\begin{equation} \begin{split} \label{23}
\sum_{n=k}^\infty T_{n,k}(x,0,0,\cdots,0|\lambda)\frac{t^n}{n!}=\frac{1}{k!}(xt)^{k}.
\end{split} \end{equation}
Thus, by \eqref{23}, we get
\begin{equation*} \begin{split} 
T_{n,k}(x,0,0,\cdots,0|\lambda)=x^{k}{0 \choose n-k}.
\end{split} \end{equation*}
From \eqref{18}, we have
\begin{equation} \begin{split} \label{24}
T_{n,k}(x,x,\cdots,x|\lambda)=x^{k}T_{n,k}(1,1,\cdots,1|\lambda),
\end{split} \end{equation}
and
\begin{equation} \begin{split} \label{25}
T_{n,k}(\alpha x_{1},\alpha x_{2},\cdots,\alpha x_{n-k+1}|\lambda)=\alpha^{k} T_{n,k}(x_{1},x_{2},\cdots,x_{n-k+1}|\lambda),
\end{split} \end{equation}
where $n,k \geq 0$ with $n-k \equiv 0$ $(mod\,\,\,2)$ and $n \geq k$.\\
Now, we observe that
\begin{equation} \begin{split} \label{26}
\textnormal{exp}&\Big(\sum_{i=1}^\infty x_{i}\big((\frac{1}{2})_{i,\lambda}-(-1)^{i}<\frac{1}{2}>_{i,\lambda}\big)\frac{t^i}{i!}\Big)\\
&=\sum_{k=0}^\infty \frac{1}{k!}\Big(\sum_{i=1}^\infty x_{i}\big((\frac{1}{2})_{i,\lambda}-(-1)^{i}<\frac{1}{2}>_{i,\lambda}\big)\frac{t^i}{i!}\Big)^{k}\\
&=1+\sum_{k=1}^\infty \frac{1}{k!}\Big(\sum_{i=1}^\infty x_{i}\big((\frac{1}{2})_{i,\lambda}-(-1)^{i}<\frac{1}{2}>_{i,\lambda}\big)\frac{t^i}{i!}\Big)^{k}\\
&=1+\sum_{k=1}^\infty \sum_{n=k}^\infty T_{n,k}(x_{1},x_{2},\cdots,x_{n-k+1}|\lambda)\frac{t^n}{n!}\\
&=1+\sum_{n=1}^\infty \big(\sum_{k=1}^{n}T_{n,k}(x_{1},x_{2},\cdots,x_{n-k+1}|\lambda)\big)\frac{t^n}{n!}.
\end{split} \end{equation}

In view of \eqref{10}, we define the degenerate central complete Bell polynomials by 
\begin{equation} \begin{split} \label{27}
\textnormal{exp}\Big(&\sum_{i=1}^\infty x_{i}\big((\frac{1}{2})_{i,\lambda}-(-1)^{i}<\frac{1}{2}>_{i,\lambda}\big)\frac{t^i}{i!}\Big)\\
&=\sum_{n=0}^\infty B_{n}^{(c)}(x_{1},x_{2},\cdots,x_{n}|\lambda)\frac{t^n}{n!}.
\end{split} \end{equation}

From \eqref{26} and \eqref{27}, we have
\begin{equation} \begin{split} \label{28}
B_{n}^{(c)}(x_{1},x_{2},\cdots,x_{n}|\lambda)=\sum_{k=0}^{n}T_{n,k}(x_{1},x_{2},\cdots,x_{n-k+1}|\lambda),
\end{split} \end{equation}

By \eqref{22} and \eqref{28}, we get
\begin{equation} \begin{split} \label{30}
B_{n}^{(c)}(x,x,\cdots,x|\lambda)&=\sum_{k=0}^{n} T_{n,k}(\underbrace{x,x,\cdots,x}_{n-k+1-times}|\lambda)\\
&=\sum_{k=0}^{n}x^{k}T_{n,k}(1,1,\cdots,1|\lambda)=\sum_{k=0}^{n} x^{k}T_{2,\lambda}(n,k)\\
&=B_{n,\lambda}^{(c)}(x).
\end{split} \end{equation}

From \eqref{26}, we note that
\begin{equation} \begin{split} \label{31}
\textnormal{exp}&\Big(\sum_{i=1}^\infty x_{i}\big((\frac{1}{2})_{i,\lambda}-(-1)^{i}<\frac{1}{2}>_{i,\lambda}\big)\frac{t^i}{i!}\Big)\\
&=1+\sum_{n=1}^\infty \frac{1}{n!}\Big(\sum_{i=1}^\infty x_{i}\big((\frac{1}{2})_{i,\lambda}-(-1)^{i}<\frac{1}{2}>_{i,\lambda}\big)\frac{t^i}{i!}\Big)^{n}\\
&=1+\frac{1}{1!}\sum_{i=1}^\infty x_{i}\big((\frac{1}{2})_{i,\lambda}-(-1)^{i}<\frac{1}{2}>_{i,\lambda}\big)\frac{t^i}{i!}\\
&+\frac{1}{2!}\Big(\sum_{i=1}^\infty x_{i}\big((\frac{1}{2})_{i,\lambda}-(-1)^{i}<\frac{1}{2}>_{i,\lambda}\big)\frac{t^i}{i!}\Big)^{2}\\
&+\frac{1}{3!}\Big(\sum_{i=1}^\infty x_{i}\big((\frac{1}{2})_{i,\lambda}-(-1)^{i}<\frac{1}{2}>_{i,\lambda}\big)\frac{t^i}{i!}\Big)^{3}+\cdots\\
&=\sum_{n=0}^\infty \bigg(\sum_{m_{1}+2m_{2}+\cdots+nm_{n}=n}\frac{n!}{m_{1}!m_{2}!\cdots m_{n}!}\\
&\times\Big(\frac{x_{1}\big((\frac{1}{2})_{1,\lambda}+<\frac{1}{2}>_{1,\lambda}\big)}{1!}\Big)^{m_{1}} \Big(\frac{x_{2}\big((\frac{1}{2})_{2,\lambda}-<\frac{1}{2}>_{2,\lambda}\big)}{2!}\Big)^{m_{2}}\times \cdots\\
&\times\Big(\frac{x_{n}\big((\frac{1}{2})_{n,\lambda}-(-1)^{n}<\frac{1}{2}>_{n,\lambda}\big)}{n!}\Big)^{m_{n}}\bigg)\frac{t^n}{n!},
\end{split} \end{equation}
where the sum is over all nonnegative integers $m_1,m_2,\cdots,m_n$ such that $m_1+2m_2+\cdots+nm_n=n$. \\ 

Now, for $n \in \mathbb{N}$ with $n \equiv 1$ $(mod\,\,\,2)$, by \eqref{27} and \eqref{31}, we get
\begin{equation} \begin{split} \label{32}
B_{n}^{(c)}&(x_{1},x_{2},\cdots,x_{n}|\lambda)=\sum_{k=1}^{n} T_{n,k}(x_{1},x_{2},\cdots,x_{n-k+1}|\lambda)\\
&=\sum_{m_{1}+2m_{2}+\cdots+nm_{n}=n}\frac{n!}{m_{1}!m_{2}!\cdots m_{n}!}\\
&\times\Big(\frac{x_{1}\big((\frac{1}{2})_{1,\lambda}+<\frac{1}{2}>_{1,\lambda}\big)}{1!}\Big)^{m_{1}} \Big(\frac{x_{2}\big((\frac{1}{2})_{2,\lambda}-<\frac{1}{2}>_{2,\lambda}\big)}{2!}\Big)^{m_{2}}\times \cdots\\
&\times\Big(\frac{x_{n}\big((\frac{1}{2})_{n,\lambda}+<\frac{1}{2}>_{n,\lambda}\big)}{n!}\Big)^{m_{n}}.
\end{split} \end{equation}

Therefore, by \eqref{32}, we obtain the following theorem.
\begin{thm} For $n \in \mathbb{N}$ with $n\equiv 1$ $(mod\,\,\,2)$, we have
\begin{equation*} \begin{split} 
B_{n}^{(c)}&(x_{1},x_{2},\cdots,x_{n}|\lambda)=\sum_{m_{1}+2m_{2}+\cdots+nm_{n}=n}\frac{n!}{m_{1}!m_{2}!\cdots m_{n}!}\\
&\times\Big(\frac{x_{1}\big((\frac{1}{2})_{1,\lambda}+<\frac{1}{2}>_{1,\lambda}\big)}{1!}\Big)^{m_{1}} \Big(\frac{x_{2}\big((\frac{1}{2})_{2,\lambda}-<\frac{1}{2}>_{2,\lambda}\big)}{2!}\Big)^{m_{2}}\times \cdots\\
&\times\Big(\frac{x_{n}\big((\frac{1}{2})_{n,\lambda}+<\frac{1}{2}>_{n,\lambda}\big)}{n!}\Big)^{m_{n}},
\end{split} \end{equation*}
where the sum is over all nonnegative integers $m_1,m_2,\cdots,m_n$ such that $m_1+2m_2+\cdots+nm_n=n$.
\end{thm}

We observe that
\begin{equation} \begin{split} \label{33}
\textnormal{exp}&\Big(x\sum_{i=1}^\infty \big((\frac{1}{2})_{i,\lambda}-(-1)^{i}<\frac{1}{2}>_{i,\lambda}\big)\frac{t^i}{i!}\Big)\\
&=1+\sum_{k=1}^\infty \frac{x^{k}}{k!}\Big(\sum_{i=1}^\infty\big((\frac{1}{2})_{i,\lambda}-(-1)^{i}<\frac{1}{2}>_{i,\lambda}\big)\frac{t^i}{i!}\Big)^{k}\\
&=1+\sum_{k=1}^\infty x^{k}\sum_{n=k}^\infty T_{n,k}(1,1,\cdots,1|\lambda)\frac{t^n}{n!}\\
&=1+\sum_{n=1}^\infty\big(\sum_{k=1}^{n} x^{k}T_{n,k}(1,1,\cdots,1|\lambda)\big)\frac{t^n}{n!}.
\end{split} \end{equation}

On the other hand,
\begin{equation} \begin{split} \label{34}
\textnormal{exp}\Big(x\sum_{i=1}^\infty \big((\frac{1}{2})_{i,\lambda}-(-1)^{i}<\frac{1}{2}>_{i,\lambda}\big)\frac{t^i}{i!}\Big)&=e^{x\big(e_{\lambda}^{\frac{1}{2}}(t)-e_{\lambda}^{-\frac{1}{2}}(t)\big)}\\
&=\sum_{n=0}^\infty B_{n,\lambda}^{(c)}(x)\frac{t^n}{n!}.
\end{split} \end{equation}

Therefore, by \eqref{33} and \eqref{34}, we obtain the following theorem.
\begin{thm} For $n,k \geq 0$ with $n \geq k$, we have
\begin{equation*} \begin{split} 
\sum_{k=0}^{n} x^{k}T_{n,k}(1,1,\cdots,1|\lambda)=B_{n,\lambda}^{(c)}(x).
\end{split} \end{equation*}
\end{thm}

By Theorem 2.5, we easily get
\begin{equation} \begin{split} \label{35}
\sum_{k=0}^{n} x^{k}T_{n,k}(1,1,\cdots,1|\lambda)=\sum_{k=0}^{n} T_{n,k}(x,x,\cdots,x|\lambda)=B_{n}^{(c)}(x,x,\cdots,x|\lambda).
\end{split} \end{equation}

\begin{cor} For $n \geq 0$, we have
\begin{equation*} \begin{split}
B_{n}^{(c)}(x,x,\cdots,x|\lambda)=B_{n,\lambda}^{(c)}(x).
\end{split} \end{equation*}
\end{cor}

\end{document}